\documentclass[11pt]{article}

\usepackage{amsmath,amsfonts,amssymb,bm,mathrsfs,amsthm,amsxtra}

\newtheorem{theorem}{Theorem}[section]
\newtheorem{remark}{Remark}[section]
\newtheorem{lemma}{Lemma}[section]
\newtheorem{proposition}{Proposition}[section]

\def\p{\partial}

\DeclareMathOperator*{\Real}{Re}
\DeclareMathOperator*{\Imag}{Im}
\DeclareMathOperator{\sgn}{sgn}

\usepackage{accents}
\makeatletter
\def\widebar{\accentset{{\cc@style\underline{\mskip10mu}}}}
\makeatother

\newtheorem{definition}{Definition}[section]

\numberwithin{equation}{section}

\begin{document}

\title{\bf Three-dimensional structural stability\\ of shock
waves in elastodynamics}

\author{{\bf Artem Shafeev}\\
Novosibirsk State University, Pirogova str. 1, 630090 Novosibirsk, Russia\\
E-mail: a.shafeev@g.nsu.ru\\[16pt]
{\bf Yuri Trakhinin}\\
Sobolev Institute of Mathematics, Koptyug av. 4, 630090 Novosibirsk, Russia\\[2pt]
and\\[2pt]
Novosibirsk State University, Pirogova str. 1, 630090 Novosibirsk, Russia\\
E-mail: trakhin@math.nsc.ru\\
}

\date{ }

{\let\newpage\relax\maketitle}

\begin{abstract}
We study the three-dimensional structural stability of shock waves for the equations of elastodynamics governing isentropic flows of compressible inviscid elastic materials. By nonlinear structural stability of a shock wave we mean the local-in-time existence and uniqueness of the discontinuous shock front solution to the hyperbolic system of elastodynamics. By using  equivalent formulations of the uniform and weak Kreiss--Lopatinski conditions for 1-shocks, we show that planar shock waves in three-dimensional elastodynamics are always at least weakly stable, and we find a condition necessary and sufficient for their uniform stability. Since the system of elastodynamics satisfies the Agranovich--Majda--Osher block structure condition, uniform stability implies structural stability of corresponding nonplanar shock waves. We also show that, as in isentropic gas dynamics, all compressive shock waves are uniformly stable for convex equations of state. This paper is a natural continuation of the previous two-dimensional analysis in \cite{MTT20,T22}. As in the two-dimensional case, we make the conclusion that the elastic force plays stabilizing role for uniform stability.
\end{abstract}



\section{Introduction}\label{sec:intro}

We consider the equations of elastodynamics governing the isentropic flow of compressible inviscid elastic materials. These equations arise as the inviscid limit of the equations of compressible viscoelasticity \cite{Daf,Gurt,Jos} of Oldroyd type \cite{Old1,Old2}. The elastodynamics equations are written as the following system of 13 conservation laws:
\begin{equation}\label{7}
\left\{
\begin{array}{l}
 \partial_t\rho  +{\rm div}\, (\rho v )=0, \\
 \partial_t(\rho v ) +{\rm div}\,(\rho v\otimes v  ) + {\nabla}p-{\rm div}\,(\rho FF^{\top})=0,\\
 \partial_t(\rho F_j ) +{\rm div}\,(\rho F_j\otimes v  -  v\otimes  \rho F_j) =0,\qquad j=1,2,3,
\end{array}
\right.
\end{equation}
where $\rho$ is the density, $v\in\mathbb{R}^3$  is the velocity, $F\in \mathbb{M}(3,3)$ is the deformation gradient,  $F_j=(F_{1j},F_{2j},F_{3j})^{\top}$ is the $j$th column of $F$, and the pressure $p=p(\rho)$ is a smooth function of $\rho$. System \eqref{7} is supplemented by the identity ${\rm div}\,(\rho F^{\top})=0$ which is the set of the three divergence constraints
\begin{equation}\label{8}
{\rm div}\,(\rho F_j)=0
\end{equation}
($j=1,2,3$) on initial data, i.e., one can show that if equations \eqref{8} are satisfied initially, then they hold for all $t > 0$. One also can show \cite{QianZhang} that the physical identity \cite{Daf}
\begin{equation}\label{8'}
\rho\det F =1
\end{equation}
is also a constraint on the initial data for the Cauchy problem. As in \cite{MTT20,T22}, we restrict ourself to the special case of neo-Hookean materials, i.e., the case when the Cauchy stress tensor has the form $\rho FF^{\top}$ corresponding to the elastic energy $W(F)=\frac{1}{2}|F|^2$ for the Hookean linear elasticity.

Using \eqref{8}, we rewrite system \eqref{7}  as
\begin{equation}
\left\{
\begin{array}{l}
{\displaystyle\frac{1}{\rho c^2}\,\frac{{\rm d} p}{{\rm d}t} +{\rm div}\,{v} =0,}\\[6pt]
{\displaystyle\rho\, \frac{{\rm d}v}{{\rm d}t}+{\nabla}p -\rho\sum_{j=1}^{3}(F_j\cdot\nabla )F_j =0 ,}\\[6pt]
\rho\,{\displaystyle \frac{{\rm d} F_j}{{\rm d}t}-\rho \,(F_j\cdot\nabla )v =0,\qquad j=1,2,3,}
\end{array}\right. \label{9}
\end{equation}
where $c^2=p'(\rho)$ is the square of the sound speed and ${\rm d} /{\rm d} t =\partial_t+({v} \cdot{\nabla} )$ is the material derivative. Equations \eqref{9} form the symmetric system
\begin{equation}
\label{10}
A_0(U )\partial_tU+\sum_{j=1}^{3}A_j(U)\partial_jU=0
\end{equation}
for $U=(p,v^\top,F_1^\top,F_2^\top,F_3^\top)^\top$, with $A_0= {\rm diag} (1/(\rho c^2) ,\rho I_{12})$ and
\[
A_j=\begin{pmatrix}
{\displaystyle\frac{v_j}{\rho c^2}} & e_j^{\top} & \underline{0}^{\top} & \underline{0}^{\top} & \underline{0}^{\top} \\[7pt]
e_j&\rho v_jI_3 & -\rho F_{j1}I_3 & -\rho F_{j2}I_3 & -\rho F_{j3}I_3 \\[3pt]
\underline{0} &-\rho F_{j1}I_3 & \rho v_jI_3 & O_3 & O_3 \\
\underline{0} &-\rho F_{j2}I_3 & O_3 & \rho v_jI_3 & O_3 \\
\underline{0} &-\rho F_{j3}I_3 & O_3 & O_3 & \rho v_jI_3
\end{pmatrix},
\]
where $\underline{0}\in\mathbb{R}^3$ is the zero vector, $e_j=(\delta_{1j},\delta_{2j},\delta_{3j})^{\top}$ are the unit vectors,  and $I_m$ and $O_m$ denote the unit and zero matrices of order $m$ respectively. In \eqref{10} we think of the density as a function of the pressure: $\rho =\rho (p)$, $c^2=1/\rho'(p)$. The symmetric system \eqref{10} is hyperbolic if   $A_0>0$, i.e.,
\begin{equation}
\rho >0,\quad  \rho '(p)>0. \label{11}
\end{equation}

Our goal in this paper is to study the structural stability of shock waves for the system of hyperbolic conservation laws \eqref{7}. By structural stability of a shock wave we mean the local-in-time existence and uniqueness of the discontinuous shock front solution to system \eqref{7}. Setting purely formally $F=0$ in \eqref{7} (strictly speaking, this contradicts to \eqref{8'}), we get the system of isentropic gas dynamics. The structural stability of shock waves in isentropic gas dynamics was proved by Majda \cite{M2} (see also \cite{Met} for a refined result) provided that the uniform stability condition introduced in \cite{M1} for the constant coefficient linearized problem with boundary conditions on a planar shock is satisfied at each point of the initial nonplanar shock wave. It should be noted that in isentropic gas dynamics there are no violently unstable planar shock waves, i.e., Hadamard-type ill-posedness examples cannot be constructed for the constant coefficient linearized problem, and shock waves are always at least weakly stable in the sense of the fulfilment of the Kreiss-Lopatinski (KL) condition \cite{Kreiss} (see Sect. \ref{s4}). The structural stability of weakly stable shock waves in isentropic gas dynamics was proved by Coulombel  and Secchi \cite{CS}.

Regarding shock waves in compressible isentropic elastodynamics (for Hookean linear elasticity), a condition sufficient for their 2D structural stability was found in \cite{MTT20} by the energy method. More precisely, by an adaptation to isentropic elastodynamics of the energy method proposed by Blokhin \cite{Bl79,Bl81,Bl82} for shock waves in full gas dynamics (see also \cite{BThand,Tsiam}), an a priori estimate without loss of derivatives for solutions to the constant coefficient linearized problem was deduced in \cite{MTT20} under a stability condition satisfied for the unperturbed rectilinear shock wave. In fact, the stability condition found in \cite{MTT20} is sufficient for uniform stability of rectilinear shock waves because the derivation of an a priori estimate without loss of derivatives for the linearized problem can be considered as an indirect test of the uniform Kreiss-Lopatinski (UKL) condition \cite{Kreiss}. With reference to \cite{M2,Met,Tsiam}, this implies the structural stability of corresponding curved shock waves in 2D isentropic elastodynamics.

At last, by the direct test of the KL and UKL conditions it was proved in \cite{T22} that the stability condition from \cite{MTT20} is not only sufficient but also necessary for the uniform stability of rectilinear shock waves in 2D isentropic elastodynamics. Moreover, it was shown in \cite{T22} that as soon as this uniform stability condition is violated, shock waves are weakly stable, i.e., as in isentropic gas dynamics, shock waves are never violently unstable. Spectral analysis in \cite{T22} is based on equivalent formulations of the KL and UKL conditions for 1-shocks proposed in \cite{Tcmp}, and its key point is a delicate study of the transition between uniform and weak stability.

For the system of hyperbolic conservation laws \eqref{7}, besides shock waves there are also characteristic discontinuities. The linear and structural stability of vortex sheets in 2D and 3D isentropic elastodynamics was studied in \cite{CHW1,CHW2,CHW3,CHW4}. We also refer to structural stability results in \cite{CHW5} for vortex sheets in nonisentropic elastodynamics (thermoelasticity). Moreover, the structural stability of contact discontinuities in thermoelasticity, for which the velocity is continuous across the discontinuity surface,  was shown in \cite{CSW} under some stability condition.

In this paper we extend the ideas of the 2D spectral analysis from \cite{T22} to the 3D case. We prove that planar shock waves in 3D isentropic elastodynamics are always at least weakly stable, and we find a condition necessary and sufficient for their uniform stability. We show that system \eqref{7} satisfies the Agranovich--Majda--Osher block structure condition \cite{Ag}. Therefore, by referring to \cite{M2}, we conclude that the uniform stability of planar elastic shocks implies structural stability of corresponding nonplanar shock waves. As in the 2D case in \cite{MTT20}, we show that all compressive shock waves are uniformly stable for convex equations of state and make the conclusion that the elastic force plays stabilizing role for uniform stability.

The plan of the rest of this paper is as follows. In Sect. \ref{s2}, we formulate the free boundary problem for shock waves in 3D elastodynamics, reduce it to that in fixed domains, discuss the Lax conditions \cite{Lax}, and prove that system \eqref{7} satisfies the Agranovich--Majda--Osher block structure condition. In Sect. \ref{s3}, we formulate the constant coefficient linearized problem associated with the nonlinear free boundary problem. In Sect. \ref{s4}, we present the main results of the paper. Sect. \ref{s5} is devoted to finding the uniform stability condition. At last,  in Sect. \ref{s6} we make some final remarks, in particular, we discuss stabilizing role of elastic force for uniform stability and show how to write down for nonplanar shocks a conterpart of the uniform stability condition obtained for a nonplanar shock wave.

\section{Free boundary problem for shock waves}
\label{s2}

Let $\Gamma (t)=\{ x_1=\varphi (t,x')\}$ be a hypersurface of strong discontinuity for system \eqref{7}, where $x'=(x_2,x_3)$. We are interested in solutions of \eqref{7} that are smooth on either side of $\Gamma (t)$. Such piecewise smooth solutions are weak solutions to the system of conservation laws \eqref{7} if and only if they satisfy the Rankine-Hugoniot jump conditions on $\Gamma (t)$:
\begin{subequations}
	\label{RH}
\begin{alignat}{2}
& [\mathfrak{m}]=0, \label{RH1}\\
& \mathfrak{m}[v_{N}] +|N|^2 [p]=\left[\rho (F_{1{N}}^2+ F_{2{N}}^2 +F_{3{N}}^2)\right],\label{RH2}\\
& \mathfrak{m}[v_{\tau}]=\left[\rho (F_{1{N}}F_{1\tau}+ F_{2{N}}F_{2\tau} +F_{3{N}}F_{3\tau})\right],\label{RH3}\\
& \mathfrak{m}[F_{j{N}}]=[\rho v_{N}F_{j{N}}],\quad j=1,2,3, \label{RH4}\\
& \mathfrak{m}[F_{j\tau}]=[\rho v_{\tau}F_{j{N}}],\label{RH5}\quad j=1,2,3,\\
& [\rho F_{j{N}}]=0,\quad j=1,2,3,\label{RH6}
\end{alignat}
\end{subequations}
where $[g]=g^+|_{\Gamma}-g^-|_{\Gamma}$ denotes the jump of $g$, with $g^{\pm}:=g$ in the domains
\[
\Omega^{\pm}(t)=\{\pm (x_1- \varphi (t,x'))>0\},
\]
and
\[
\mathfrak{m}^{\pm}=\rho^{\pm} (v_{N}^{\pm}-\partial_t\varphi),\quad v_{N}^{\pm}=
v^\pm\cdot N,\quad N=(1,-\p_2\varphi,-\p_3\varphi)^{\top}, \quad
F_{j{N}}^{\pm}=F_{j}^{\pm}\cdot N ,
\]
\[
v_{\tau}^{\pm}=(v_1^{\pm}\partial_2\varphi +v_2^{\pm},v_1^{\pm}\partial_3\varphi +v_3^{\pm})^{\top}, \quad
F_{j\tau}^{\pm}=(F_{1j}^{\pm}\partial_2\varphi +F_{2j}^{\pm},F_{1j}^{\pm}\partial_3\varphi +F_{3j}^{\pm})^{\top},
\]
and $\mathfrak{m}:=\mathfrak{m}^{\pm}|_{\Gamma}$ is the mass transfer flux across the discontinuity. The jump conditions \eqref{RH6} come actually from the divergence constraints \eqref{8}. On the other hand, conditions \eqref{RH4} are rewritten as $\partial_t\varphi [\rho F_{j{\rm N}}]=0$. That is, conditions \eqref{RH4} are implied by \eqref{RH6} and can be thus excluded from system \eqref{RH}.

Mathematically, {\it shock waves} are noncharacteristic discontinuities, i.e., for them the hypersurface $\Gamma (t)$ is a not a characteristic of the hyperbolic system \eqref{10}. As in gas dynamics, one can show that this is so if $\mathfrak{m}\neq 0$ and $[\rho ]\neq 0$. Simple arguments for deriving the boundary conditions for 3D shock waves from \eqref{RH} are essentially the same as in the 2D case in \cite{MTT20}, but for the reader's convenience we present them below.

In view of \eqref{RH6}, conditions \eqref{RH3} and \eqref{RH5} form the following linear algebraic system for the jumps $[v_{\tau}]$, $[F_{1\tau}]$, $[F_{2\tau}]$ and $[F_{3\tau}]$:
\begin{equation}\label{m0}
\begin{pmatrix}
\mathfrak{m} & -\rho^+ F_{1{N}}^+ & -\rho^+ F_{2{N}}^+ & -\rho^+ F_{3{N}}^+\\
-\rho^+ F_{1{N}}^+ & \mathfrak{m} & 0 & 0\\
-\rho^+ F_{2{N}}^+ & 0 & \mathfrak{m} & 0\\
-\rho^+ F_{3{N}}^+ & 0 &  0 & \mathfrak{m}
\end{pmatrix}
\begin{pmatrix} [v_{\tau}]\\ [F_{1\tau}]\\ [F_{2\tau}]\\ [F_{3\tau}]\end{pmatrix} =0\quad\mbox{on}\ \Gamma .
\end{equation}
Since $\mathfrak{m}\neq 0$, this system has a nonzero solution if
\[
\mathfrak{m}^2=(\rho^+)^2\left.\left((F_{1{N}}^+)^2+(F_{2{N}}^+)^2 +(F_{3{N}}^+)^2 \right)\right|_{\Gamma}.
\]
We assume that $\mathfrak{m}^2\neq(\rho^+)^2\left.\left((F_{1{N}}^+)^2+(F_{2{N}}^+)^2 +(F_{3{N}}^+)^2 \right)\right|_{\Gamma}\,$, i.e.,
\begin{equation}\label{m}
(v_{N}^+-\partial_t\varphi)^2\neq (F_{1{N}}^+)^2+(F_{2{N}}^+)^2 +(F_{3{N}}^+)^2 \quad\mbox{on}\ \Gamma .
\end{equation}
We will see below that \eqref{m} holds thanks to the Lax conditions \cite{Lax}.

From \eqref{m0} and \eqref{m} we obtain
\[
[v_{\tau}]=0,\quad  [F_{j\tau}]=0,\quad j=1,2,3.
\]
In view of  \eqref{RH1} and \eqref{RH6}, we rewrite condition \eqref{RH2} as
\[
\mathfrak{M}[V]+|N|^2[p]=0,
\]
where
\[
\mathfrak{M}=\mathfrak{m}^2-(\rho^+)^2\left.\left((F_{1{N}}^+)^2+(F_{2{N}}^+)^2 +(F_{3{N}}^+)^2 \right)\right|_{\Gamma},\quad V^\pm=1/\rho^\pm
\]
($\mathfrak{M}\neq 0$, cf. \eqref{m}). We thus obtain the following 13 boundary conditions on $\Gamma (t)$:
\begin{equation}\label{bc}
\begin{split}
& [\mathfrak{m}]=0,\quad \mathfrak{M}[V]+|N|^2[p]=0,\quad [v_{\tau}]=0,\\
&  [F_{j\tau}]=0,\quad [\rho F_{j{N}}]=0,\quad j=1,2,3.
\end{split}
\end{equation}

The free boundary problem for shock waves is the problem for the systems
\begin{equation}
A_0(U^{\pm})\partial_tU^{\pm}+\sum_{j=1}^{3}A_j(U^{\pm} )\partial_jU^{\pm}=0\quad \mbox{in}\ \Omega^{\pm}(t)
\label{21}
\end{equation}
with the boundary conditions \eqref{bc} on $\Gamma (t)$ and  the initial data
\begin{equation}
{U}^{\pm} (0,{x})={U}_0^{\pm}({x}),\quad {x}\in \Omega^{\pm} (0),\quad \varphi (0,x')=\varphi _0(x'),\quad x'\in\mathbb{R}^{2}.\label{indat}
\end{equation}
The initial data \eqref{indat} should satisfy not only the hyperbolicity conditions \eqref{11} but also constraints \eqref{8} and \eqref{8'}. As for the Cauchy problem, we can show that these constraints are preserved by problem \eqref{bc}--\eqref{indat} (the proof of this is interely the same as for the 2D case in \cite{MTT20}).
The {\it structural stability} of a shock wave means the local-in-time existence and uniqueness of a smooth solution $(U^+,U^-,\varphi)$ to  problem \eqref{bc}--\eqref{indat}.

The free boundary problem  \eqref{bc}--\eqref{indat} is reduced to that in the fixed domains $[0,T]\times\mathbb{R}^3_{\pm}$ by the simple change of variables ${x}^\sharp_1=x_1-\varphi (t,x')$, where $\mathbb{R}^3_{\pm}=\{\pm x_1>0,\ x'\in\mathbb{R}^{2}\}$. Dropping the superscript $``\sharp"$, from systems \eqref{21} we obtain
\begin{equation}
A_0(U^{\pm})\partial_tU^{\pm}+\widetilde{A}_{1}(U^{\pm},\varphi )\partial_1U^{\pm}+A_2(U^{\pm} )\partial_2U^{\pm}+A_3(U^{\pm} )\partial_3U^{\pm}=0\qquad \mbox{in }[0,T]\times\mathbb{R}^3_{\pm},
\label{23}
\end{equation}
where
\begin{equation}
\widetilde{A}_{1}=\widetilde{A}_{1}(U,\varphi )= A_1(U)-A_0(U)\partial_t\varphi -A_2(U)\partial_2\varphi -A_3(U)\partial_3\varphi .
\label{AN}
\end{equation}
The boundary conditions for \eqref{23} are \eqref{bc} on the plain $x_1=0$.

As is known, the Lax conditions \cite{Lax} guarantee the correct number of boundary conditions on $x_1=0$ for systems \eqref{23}. The Lax conditions are the inequalities
\[
\lambda_{k-1}^- <\partial_t\varphi <\lambda_k^-,\quad \lambda_k^+ <\partial_t\varphi <\lambda_{k+1}^+
\]
which must hold for some fixed integer $k$, where in our case $1\leq k\leq  =13$ and $\lambda_j^\pm$ ($j=\overline{1,{13}}$) are the eigenvalues of the matrices
\begin{equation}
A_{N}^{\pm}:=\left(A_0(U^{\pm})\right)^{-1}\left.\left( A_1(U^{\pm})-A_2(U^{\pm})\partial_2\varphi -A_3(U^{\pm})\partial_3\varphi \right)\right|_{x_1=0},
\label{bmat}
\end{equation}
with $U^{\pm}$ and $\varphi$ satisfying the boundary conditions \eqref{bc} on $x_1=0$. Moreover, $\lambda_j^\pm$ are numbered as
\[
\lambda_1^-\leq \ldots \leq\lambda_{13}^-,\quad \lambda_1^+\leq \ldots \leq\lambda_{13}^+,
\]
and we take $\lambda^-_0:=-2|\partial_t\varphi |$ and $\lambda_{14}^+:=2|\partial_t\varphi |$. Direct computations give
\begin{equation}
\left\{
\begin{array}{l}
 \lambda_1^{\pm}=v_{N}^{\pm}-\sqrt{(\bar{c}^{\pm})^2+(F_{1{N}}^{\pm})^2+(F_{2{N}}^{\pm})^2 +(F_{3{N}}^{\pm})^2}\,,\\[6pt] \lambda_2^{\pm}=\lambda_3^{\pm}=v_{N}^{\pm}-\sqrt{(F_{1{N}}^{\pm})^2+(F_{2{N}}^{\pm})^2
 +(F_{3{N}}^{\pm})^2}\,,\\[6pt]  \lambda_{4}^{\pm}=\ldots=\lambda_{10}^{\pm}=v_{N}^{\pm}\,,\\[6pt]
 \lambda_{11}^{\pm}= \lambda_{12}^{\pm}=v_{N}^{\pm}+\sqrt{(F_{1{N}}^{\pm})^2+(F_{2{N}}^{\pm})^2
 +(F_{3{N}}^{\pm})^2}\,,\\[6pt] \lambda_{13}^{\pm}=v_{N}^{\pm}+\sqrt{(\bar{c}^{\pm})^2+(F_{1{N}}^{\pm})^2+
 (F_{2{N}}^{\pm})^2 +(F_{3{N}}^{\pm})^2}\qquad \mbox{on}\ x_1=0,
\end{array}
\right.
\label{24}
\end{equation}
where $\bar{c}^{\pm}= c^{\pm}|N|$ and $c^{\pm} =1/\sqrt{\rho'(p^{\pm})}$ are the sound speeds ahead and behind of the shock.

Without loss of generality we assume that $v_{N}^{\pm}|_{x_1=0}>\partial_t\varphi$. Then, we can check that, as in gas dynamics (as well as in 2D elastodynamics \cite{MTT20,T22}), only 1-shocks are possible:
\begin{equation}
\begin{split}
& v_{N}^--\partial_t\varphi >\sqrt{(\bar{c}^-)^2+(F_{1{N}}^{-})^2+(F_{2{N}}^{-})^2 +(F_{3{N}}^{-})^2}\,, \\[6pt]
& \sqrt{(F_{1{N}}^{+})^2+(F_{2{N}}^{+})^2+(F_{3{N}}^{+})^2}   <v_N^+-\partial_t\varphi
 \\[6pt]
& \qquad\qquad\qquad\qquad\qquad\qquad\ \,
 <\sqrt{(\bar{c}^{+})^2+(F_{1{N}}^{+})^2+(F_{2{N}}^{+})^2 +(F_{3{N}}^{+})^2}  \qquad \mbox{on}\ x_1=0.
\end{split}
\label{25}
\end{equation}
We see that the Lax conditions \eqref{25} imply the fulfilment of
the above assumption \eqref{m}.

\begin{remark}
\label{r1}
{\rm The boundary matrices in \eqref{bmat} are nothing else than
$\mathcal{A}(U^\pm,N )|_{x_1=0}$, where
\[
\mathcal{A}(U,\xi )= \left(A_0(U)\right)^{-1}\sum_{j=1}^3\xi_jA_j(U),
\]
with $\xi \in\mathbb{R}^3\setminus\{\underline{0}\}$. We can easily calculate the eigenvalues $\lambda_j=\lambda_j(U,\xi)$ of $\mathcal{A}(U,\xi )$, cf. \eqref{24}:
\[
\begin{array}{l}
 \lambda_1=v\cdot\xi-\sqrt{{c}^2|\xi |^2+(F_1\cdot\xi)^2+(F_2\cdot\xi)^2+(F_3\cdot\xi)^2}\,,\\[6pt] \lambda_2=\lambda_3=v\cdot\xi-\sqrt{(F_1\cdot\xi)^2+(F_2\cdot\xi)^2+(F_3\cdot\xi)^2}\,,\\[6pt]  \lambda_{4}=\ldots=\lambda_{10}=v\cdot\xi,\\[6pt]
 \lambda_{11}= \lambda_{12}=v\cdot\xi+\sqrt{(F_1\cdot\xi)^2+(F_2\cdot\xi)^2+(F_3\cdot\xi)^2}\,,\\[6pt] \lambda_{13}=v\cdot\xi+\sqrt{{c}^2|\xi |^2+(F_1\cdot\xi)^2+(F_2\cdot\xi)^2+(F_3\cdot\xi)^2}\,.
\end{array}
\]
Since $\det F>0$ (see \eqref{8'}, \eqref{11}),
\[
(F_1\cdot\xi)^2+(F_2\cdot\xi)^2+(F_3\cdot\xi)^2> 0\quad \mbox{for all}\quad\xi \in\mathbb{R}^3\setminus\{\underline{0}\}.
\]
Therefore, the eigenvalues $\lambda_j(U,\xi)$ are of constant multiplicity for all $\xi \in\mathbb{R}^3\setminus\{\underline{0}\}$. Moreover, since system \eqref{10} is symmetric hyperbolic, the matrix $\mathcal{A}(U,\xi )$ is diagonalizable, i.e., all the multiple eigenvalues $\lambda_j(U,\xi)$ are semi-simple. Referring to \cite{MZ}, we may conclude that system \eqref{10} satisfies the Agranovich--Majda--Osher {\it block structure condition} \cite{Ag,M1,MajOsh}. In view of the results in \cite{M1,M2,Met}, this means that the uniform linearized stability of a planar shock wave in elastodynamics  implies the nonlinear structural stability of corresponding nonplanar shock waves.
}
\end{remark}

\section{Constant coefficient linearized problem associated with planar shock waves}
\label{s3}

Without loss of generality we assume that the unperturbed planar shock wave is given by the equation $x_1=0$. We now consider a constant solution $(U^+,U^-,\varphi) =(\widehat{U}^+,\widehat{U}^-,0)$ of systems \eqref{23} and the boundary conditions \eqref{bc} associated with the planar shock wave $x_1=0$:
\[
\widehat{U}^{\pm}=(\hat{p}^{\pm},\hat{v}^{\pm},\widehat{F}_1^{\pm},\widehat{F}_2^{\pm},
\widehat{F}_3^{\pm}), \quad \hat{\rho}^{\pm}=\rho (\hat{p}^{\pm}) >0,\quad \hat{c}^{\pm}=1/\sqrt{\rho'(\hat{p}^{\pm})}>0,
\]
\[
\hat{v}^{\pm}=(\hat{v}_1^{\pm},\hat{v}_2^{\pm} ,\hat{v}_3^{\pm}),\quad
\widehat{F}_j^{\pm} =(\widehat{F}_{1j}^{\pm},\widehat{F}_{2j}^{\pm} ,\widehat{F}_{3j}^{\pm}), \quad j=1,2, 3,
\]
where all the hat values are given constants.  In view of the third condition in \eqref{bc}, $\hat{v}_l^+=\hat{v}_l^-$  for $l=2,3$ and we can choose a reference frame in which
$\hat{v}_2^\pm=\hat{v}_3^\pm=0$. The rest constants satisfy the relations
\begin{equation}
\begin{split}
& \frac{\hat{\rho}^+}{\hat{\rho}^-}=\frac{\hat{v}_1^-}{\hat{v}_1^+},\quad \frac{\hat{\rho}^+}{\hat{\rho}^-}\big\{ (\hat{v}_1^+)^2-\big((\widehat{F}_{11}^+)^2+
(\widehat{F}_{12}^+)^2 +(\widehat{F}_{13}^+)^2\big)\big\}[\hat{\rho}]=[\hat{p}],\\ & \big[\widehat{F}_{lj}\big]=0,\quad l=2,3,\quad \big[\hat{\rho}\widehat{F}_{1j}\big]=0
\end{split}
\label{sbc}
\end{equation}
following from \eqref{bc}, where $j=1,2,3$, $[\hat{\rho}]=\hat{\rho}^+-\hat{\rho}^-$, etc., and without loss of generality we assume that $\hat{v}_1^{\pm}>0$.

The Lax conditions \eqref{25} for the constant solution read:
\begin{equation}
M_->\frac{M}{\sqrt{M^2-M_1^2}},
\label{Mach-}
\end{equation}
\begin{equation}
M_1<M<M_*,
\label{Mach}
\end{equation}
where $M_-={\hat{v}_1^-}/{\hat{c}^-}$ and $M={\hat{v}_1^+}/{\hat{c}^+}$ are the upstream and downstream  Mach numbers respectively,
\[
M_1=|\widebar{\mathcal{F}}_{1}|,\quad M_*=\sqrt{1+|\widebar{\mathcal{F}}_{1}|^2},\quad
\mathcal{F}_{ij}=\widehat{F}_{ij}^+/\hat{c}^+,\quad i,j=1,2,3,
\]
$\mathcal{F}_{ij}$ are the components of the unperturbed scaled deformation gradient $\mathcal{F}=(\mathcal{F}_{ij})_{i,j=1,2,3}$ behind of the shock, $\mathcal{F}_j=(\mathcal{F}_{1j},\mathcal{F}_{2j},\mathcal{F}_{3j})^\top$ and $\widebar{\mathcal{F}}_j=({\mathcal{F}}_{j1},{\mathcal{F}}_{j2},
{\mathcal{F}}_{j3})^\top$ ($j=1,2,3$) are the columns and (transposed) rows of the matrix $\mathcal{F}$ respectively.
Inequality \eqref{Mach-} follows from the first inequality in \eqref{25} and relations \eqref{sbc}.

For 1-shocks the linearized system for the perturbation $\delta U^-$ ahead of the planar shock (for $x_1<0$) does not need boundary conditions because all the characteristics of  are incoming in the shock. Therefore, without loss of generality we may assume that $\delta U^-\equiv 0$. Following \cite{MTT20,T22}, we write down the constant coefficient linearized  problem in a dimensionless form for the scaled perturbation $U=(p,v^\top,F_1^\top,F_2^\top,F_3^\top)^\top$ behind of the shock wave resulting from the linearization of \eqref{bc}, \eqref{23} about the constant solution described above:
\begin{equation}
\left\{
\begin{array}{l}
{\displaystyle
Lp+{\rm div}\,v=0,} \\[3pt]
{\displaystyle M^2Lv+\nabla p -\sum_{j=1}^3(\mathcal{F}_j\cdot\nabla )F_j=0, }\\[12pt]
{\displaystyle LF_j- (\mathcal{F}_j\cdot\nabla )v=0,\quad j=1,2,3, \qquad \mbox{for}\ x_1>0,}
\end{array}\right. \label{ls}
\end{equation}
\begin{equation}
\left\{
\begin{array}{l}
{\displaystyle v_1 +d_0p -\frac{1}{M^2R}(\ell_2v_2+\ell_3v_3)=0,  }\\[12pt]
a_0p+(1-R)\partial_{\star}\varphi =0,\qquad
v_k +(1-R)\partial_k\varphi= 0, \\[6pt]
{\displaystyle  F_{1j}+\mathcal{F}_{1j}\,p-\frac{1}{R}(\mathcal{F}_{2j}v_2+
\mathcal{F}_{3j}v_3)=0,\qquad F_{kj}-
\mathcal{F}_{1j}\,v_k=0,
  \qquad \mbox{on}\ x_1=0, }\\[12pt]
 j=1,2,3,\quad k=2,3,
\end{array}\right. \label{lbc}
\end{equation}
\begin{equation}
{U} (0,{x})={U}_0({x}),\quad {x}\in \mathbb{R}^3,\quad \varphi (0,{x}')=\varphi _0({x}'),\quad {x}'\in\mathbb{R}^{2},\label{lindat}
\end{equation}
where
\[
L =\partial_t+\partial_1,\quad
d_0=\frac{M_*^2+M^2}{2M^2},\quad R=\frac{\hat{\rho}^+}{\hat{\rho}^-},\quad
\ell_2=\widebar{\mathcal{F}}_{1}\cdot\widebar{\mathcal{F}}_{2}\,,\quad
\ell_3=\widebar{\mathcal{F}}_{1}\cdot\widebar{\mathcal{F}}_{3}\,,
\]
\[
a_0=-\frac{\beta^2R}{2M^2},\quad
\beta =\sqrt{M_*^2-M^2},\quad
\partial_{\star}=\partial_t-\frac{1}{M^2}(\ell_2\partial_2 +\ell_3\partial_3),
\]
and the following scaled values
\[
\tilde{x}=\frac{x}{\hat{l}},\quad \tilde{t}=\frac{\hat{v}_1^+t}{\hat{l}},\quad p=\frac{\delta p^+}{\hat{\rho}^+(\hat{c}^+)^2},\quad v=\frac{\delta v^+}{\hat{v}_1^+},\quad F_{ij}=\frac{\delta F_{ij}^+}{\hat{c}^+},\quad  \varphi =\frac{\delta\varphi}{\hat{l}},
\]
are used instead of the original (unscaled) perturbation $\delta U^+=(\delta p^+,\delta v^+,\delta F_1^+,\delta F_2^+,\delta F_3^+)$ behind the shock and the shock perturbation $\delta\varphi$, with $\hat{l}$ being a typical length (the tildes were dropped in \eqref{ls}--\eqref{lindat}).

System \eqref{ls}  can rewritten  as
\begin{equation}\label{lsm}
\mathcal{A}_0\partial_tU+\sum_{j=1}^3\mathcal{A}_j\partial_jU=0\quad\mbox{for}\ x\in \mathbb{R}^3_+,
\end{equation}
where $\mathcal{A}_0={\rm diag}\,(1,M^2I_3,I_{9})$ and
\[
\mathcal{A}_j=\delta_{1j}\mathcal{A}_0 +\begin{pmatrix}
0 & e_j^{\top} & \underline{0}^{\top} & \underline{0}^{\top} &\underline{0}^{\top}  \\[7pt]
e_j&O_3 &  -\mathcal{F}_{j1}I_3 & -\mathcal{F}_{j2}I_3 & - \mathcal{F}_{j3}I_3  \\[3pt]
\underline{0} &- \mathcal{F}_{j1}I_3 & O_3 & O_3 & O_3 \\
\underline{0} &- \mathcal{F}_{j2}I_3 & O_3 & O_3 & O_3 \\
\underline{0} &- \mathcal{F}_{j3}I_3 & O_3 & O_3 & O_3
\end{pmatrix}.
\]
By cross differentiation we exclude $\varphi$ from the boundary conditions \eqref{lbc}:
\begin{equation}\label{lbc2'}
\partial_{\star}v_k=a_0\partial_kp,\quad k=2,3\qquad \mbox{on}\ x_1=0,
\end{equation}
and can thus rewrite the boundary conditions in  the form
\begin{equation}
\mathfrak{B}_0\partial_tU+\mathfrak{B}_2\partial_2U+\mathfrak{B}_3\partial_3U+
\mathfrak{B}_{1}U =0\quad \mbox{on}\ x_1=0, \label{72}
\end{equation}
where the matrices $\mathfrak{B}_{\alpha}$ ($\alpha =\overline{0,3}$) of order $12\times 13$ can be easily written down.

\section{Main results}
\label{s4}

The KL and UKL conditions for initial boundary value problems for linear constant coefficient hyperbolic systems were first introduced by Kreiss \cite{Kreiss}. In view of the Lax conditions \eqref{Mach}, problem \eqref{lindat}, \eqref{lsm},  \eqref{72} has the 1-shock property in the sense the nonsingular matrix $(\mathcal{A}_0)^{-1}\mathcal{A}_1$ has a single ne\-gative eigenvalue. For such kind of noncharacteristic initial boundary value problems for symmetric hyperbolic systems equivalent formulations of the KL and UKL conditions were given in \cite{Tcmp} by using ideas of Gardner and Kruskal \cite{GK}. Roughly speaking, the main idea of \cite{Tcmp} borrowed from \cite{GK} is to carry all the calculations over the single characteristic incoming in the shock in the region behind of the shock wave.

Following \cite{Tcmp} (see also \cite[Appendix A]{T22}), we write down the following linear algebraic system for a vector $X$ associated with problem \eqref{lsm}, \eqref{72}:
\begin{align}
& (s\mathcal{A}_0+\lambda \mathcal{A}_1+{i}\omega_2 \mathcal{A}_2
+{i}\omega_3 \mathcal{A}_3)X =0\,,\label{3.5}
\\
& \big(\mathcal{A}_1\widetilde{U}_0\big)\cdot X=0,
\label{3.6}
\end{align}
where
\[
s=\eta +{i}\xi ,\quad \eta >0,\quad (\xi ,{\omega}')\in \mathbb{R}^3,\quad \omega'=(\omega_2,\omega_3)
\]
($s$ and $\omega'$ are the Laplace and Fourier variables respectively), the vector $\widetilde{U}_0=\widetilde{U}_0(\eta ,\xi ,\omega' )$ satisfies the Fourier--Laplace transform of the boundary conditions \eqref{72},
\begin{equation}
(s\mathfrak{B}_0+i\omega_2\mathfrak{B}_2+i\omega_3\mathfrak{B}_3+\mathfrak{B}_1)\widetilde{ U}_0=0,\label{3.2}
\end{equation}
and $\lambda=\lambda^+ (\eta ,\xi ,\omega' )$ is the single root $\lambda$ of the dispersion relation
\begin{equation}
\det (s\mathcal{A}_0+\lambda \mathcal{A}_1+{i}\omega_2 \mathcal{A}_2+{i}\omega_3 \mathcal{A}_3)=0 \label{3.4}
\end{equation}
lying strictly in the open right-half complex plane ($\Real\lambda^+ >0$).

Since $\lambda^+$ is a simple root, we can choose 12 linearly independent equations of system \eqref{3.5}. Adding them to equation \eqref{3.6}, we obtain for the vector $X$ the linear algebraic system
\begin{equation}
\mathfrak{L}X =0
\label{X}
\end{equation}
whose determinant is, in fact, the Lopatinski determinant (see \cite{Tcmp} and \cite[Appendix A]{T22}). Following \cite{Tcmp,T22}, we are now ready to give the definitions of the KL and UKL conditions for problem \eqref{lsm}, \eqref{72}.

\begin{definition}
Problem \eqref{lsm}, \eqref{72} satisfies the KL condition if $\det \mathfrak{L} (\eta
,\xi ,\omega' ,\lambda^+ ) \neq 0$ for all $\eta >0$ and  $(\xi , \omega' )\in
\mathbb{R}^3$.
\label{d3.2}
\end{definition}

\begin{definition}
Problem \eqref{lsm}, \eqref{72} satisfies the UKL condition if $\det \mathfrak{L} (\eta
,\xi ,\omega' ,\lambda^+ ) \neq 0$ for all $\eta \geq 0$ and  $(\xi , \omega' )\in
\mathbb{R}^3$ (with $\eta^2+\xi^2+|\omega'|^2\neq 0$), where $\lambda^+ (0,\xi ,\omega' )= \lim\limits_{\eta \rightarrow +0}
\lambda^+ (\eta ,\xi ,\omega' )  $.
\label{d3.3}
\end{definition}

The planar shock wave for which problem \eqref{lsm}, \eqref{72} satisfies the UKL condition (under suitable assumptions on the constant states $\widehat{U}^+$ and $\widehat{U}^-$) is called {\it uniformly stable}. If the KL condition holds, but the UKL is violated, then the planar shock is called {\it weakly stable}. If even the KL condition is violated, then the shock wave is {\it violently unstable} meaning that the linearized stability problem \eqref{lsm}, \eqref{72} is ill-posed.

We are now in a position to formulate the main result of this paper.

\begin{theorem}
Let the planar shock wave $x_1=0$ in compressible isentropic elastodynamics satisfies the Lax conditions \eqref{Mach-} and \eqref{Mach}. Then, this shock wave is uniformly stable if and only if
\begin{multline}
\min_{\bar{\omega}_2^2+\bar{\omega}_3^2=1}
\bigg(M  \sqrt{1+|\widebar{\mathcal{F}}_1|^2+
|\bar{\omega}_2\widebar{\mathcal{F}}_2+\bar{\omega}_3\widebar{\mathcal{F}}_3|^2+
\big|\widebar{\mathcal{F}}_1\times (\bar{\omega}_2\widebar{\mathcal{F}}_2+\bar{\omega}_3\widebar{\mathcal{F}}_3)\big|^2}
\\
 -\big|\widebar{\mathcal{F}}_1\cdot(\bar{\omega}_2\widebar{\mathcal{F}}_2+
 \bar{\omega}_3\widebar{\mathcal{F}}_3)\big|\sqrt{1+|\widebar{\mathcal{F}}_1|^2
-M^2}\,\bigg)^2 -
\big( 1 + |\widebar{\mathcal{F}}_1|^2\big)^2 \big(R (M^2 -  |\widebar{\mathcal{F}}_1|^2) -1 \big) >0,
\label{usc'}
\end{multline}
where $M$ is the downstream Mach number, $R$ measures the competition between downstream and upstream densities and $\widebar{\mathcal{F}}_j$ ($j=1,2,3$) are the transposed rows of the scaled deformation gradient $\mathcal{F}=(\mathcal{F}_{ij})_{i,j=1,2,3}$ behind of the shock. If \eqref{usc'} is violated, then the planar shock wave is weakly stable.
\label{t1}
\end{theorem}

As was noted above (see Remark \ref{r1}), the uniform stability of a planar shock wave in elastodynamics  implies the nonlinear structural stability of corresponding nonplanar shock waves, i.e., uniform stability is sufficient for structural stability. For completeness, we now formulate the corresponding sufficient structural stability condition which should be satisfied by the initial data of the free boundary problem \eqref{bc}--\eqref{indat} at each point of the initial nonplanar Lax shock wave $\Gamma (0)$.

\begin{theorem}
Let the initial shock wave with the equation $x_1=\varphi_0 (x')$ be a Lax shock, i.e., the initial data \eqref{indat} of the free boundary problem \eqref{bc}--\eqref{indat} satisfy the inequalities in \eqref{25} on $\Gamma (0)$. Let also the initial data satisfy the hyperbolicity conditions \eqref{11} on either side of $\Gamma (0)$ as well as suitable compatibility conditions. Then the condition
\begin{multline}
\min_{\bar{\omega}_2^2+\bar{\omega}_3^2=1}
\bigg( \mathcal{M}  \sqrt{1+|\widebar{F}_N|^2+
|\bar{\omega}_2\widebar{F}_2+\bar{\omega}_3\widebar{F}_3|^2+
\big|\widebar{F}_N\times (\bar{\omega}_2\widebar{F}_2+\bar{\omega}_3\widebar{F}_3)\big|^2}
\\
-\big|\widebar{F}_N\cdot(\bar{\omega}_2\widebar{F}_2+
 \bar{\omega}_3\widebar{F}_3)\big|\sqrt{1+|\widebar{F}_N|^2
-\mathcal{M}^2}\,\bigg)^2 \\ -
\big( 1 + |\widebar{F}_N|^2\big)^2 \big(\mathcal{R} (\mathcal{M}^2 -  |\widebar{F}_N|^2) -1 \big)   \geq \varepsilon >0
\label{strc}
\end{multline}
satisfied by the initial data \eqref{indat} at each point of the hypersurface $x_1=\varphi_0 (x')$ is sufficient for structural stability, i.e., for the existence of a unique local-in-time smooth solution $(U^+,U^-,\varphi)$ to  problem \eqref{bc}--\eqref{indat}, where
\begin{equation}
\mathcal{M}= \frac{v_N^+-\p_t\varphi}{{c}^+|N|},\quad \widebar{F}_N =\frac{1}{{c}^+|N|}\,(F_{1N}^+,
F_{2N}^+,F_{3N}^+)^{\top},\quad \widebar{F}_k=\frac{1}{c^+}\,(F^+_{k1},F^+_{k2},F^+_{k3})^{\top},\quad k=2,3,
\label{ff}
\end{equation}
and $\mathcal{R}=\rho^+/\rho^-$.
\label{t2}
\end{theorem}

\begin{remark}
\label{r2}
{\rm
The uniform stability condition found in \cite{MTT20,T22} for the 2D case can be rewritten in the form, which is in some sense similar to \eqref{usc'}:
\begin{multline}
\bigg(M  \sqrt{1+|\widebar{\mathcal{F}}_1|^2+
|\widebar{\mathcal{F}}_2|^2+(\det \mathcal{F})^2}
 -\big|\widebar{\mathcal{F}}_1\cdot\widebar{\mathcal{F}}_2\big|\sqrt{1+|\widebar{\mathcal{F}}_1|^2
-M^2}\,\bigg)^2  \\ -
\big( 1 + |\widebar{\mathcal{F}}_1|^2\big)^2 \big(R (M^2 -  |\widebar{\mathcal{F}}_1|^2) -1 \big) >0,
\label{usc2D}
\end{multline}
where we use natural 2D counterparts of notations exploited above in the 3D case. It can be easily seen that by setting formally $\mathcal{F}=0$ in  \eqref{usc'} and \eqref{usc2D} both conditions are reduced to the uniform stability condition
\begin{equation}
M^2(R-1)<1
\label{gas}
\end{equation}
found by Majda \cite{M1} (and written in our notations) for shock waves in isentropic gas dynamics. Unfortunately, we cannot totally exclude the ``wave vector'' $\bar{\omega}'=(\bar{\omega}_2,\bar{\omega}_3)$ from the 3D condition \eqref{usc'}, but one can do this for particular deformations, e.g., for the case of pure {\it stretching} $\mathcal{F}={\rm diag}( \mathcal{F}_{11},\mathcal{F}_{22},\mathcal{F}_{33})$. The 2D  uniform stability condition \eqref{usc2D} for the case of stretching $\mathcal{F}={\rm diag}( \mathcal{F}_{11},\mathcal{F}_{22})$ is reduced to
\begin{equation}\label{str}
1+\mathcal{F}_{11}^2 +M^2 -R(1+\mathcal{F}_{11}^2)(M^2-\mathcal{F}_{11}^2)+\mathcal{F}_{22}^2M^2>0.
\end{equation}
It follows from \eqref{usc'} that the 3D version of \eqref{str} reads
\[
1+\mathcal{F}_{11}^2 +M^2 -R(1+\mathcal{F}_{11}^2)(M^2-\mathcal{F}_{11}^2)+
\min (\mathcal{F}_{22}^2,\mathcal{F}_{33}^2)\,M^2>0.
\]

}
\end{remark}

\section{Spectral analysis: proof of Theorem 4.1}
\label{s5}

We first compute the dispersion relation \eqref{3.4} for finding the single root $\lambda=\lambda^+ (\eta ,\xi ,\omega' )$ (with $\Real\lambda^+ >0$):
\begin{equation}
\Omega^7\left( M^2\Omega^2-\sigma_1^2-\sigma_2^2-\sigma_3^2\right)^2\left( M^2\Omega^2-\sigma_1^2-\sigma_2^2-\sigma_3^2-\lambda^2+\omega^2\right)=0,
\label{3.4'}
\end{equation}
where $\Omega =s+\lambda$, $\omega =|\omega'|$, $\sigma_j=\mathcal{F}_j\cdot\zeta$,  and $\zeta =(\lambda,i\omega')$, with $j=1,2,3$. It is clear that $\Omega \neq 0$ for $\lambda=\lambda^+$. Moreover,
\begin{equation}
 M^2\Omega^2-\sigma_1^2-\sigma_2^2-\sigma_3^2\neq 0\quad\mbox{for }\lambda=\lambda^+.
 \label{neq0}
\end{equation}
Indeed, for $\omega' =0$ the last expression vanishes at $\lambda=-Ms/(M\pm M_1)$ for which $\Real\lambda<0$. Referring to Hersh's lemma \cite{Hersh}, we conclude that \eqref{neq0} is true also for all $\omega'\neq 0$ (and $\eta >0$). Hence, it follows from  \eqref{3.4'} that $\lambda^+$ is one of the two roots of the equation
\[
 M^2\Omega^2-\sigma_1^2-\sigma_2^2-\sigma_3^2-\lambda^2+\omega^2=0,
\]
which can be rewritten as
\begin{equation}\label{83'}
M^2\Omega^2-M_*^2\lambda^2 +K_0\omega^2=2i\ell_0\lambda\omega ,
\end{equation}
where
\begin{align}
& K_0=1+M_2^2,\quad M_2=|\bar{\omega}_2\bar{\mathcal{F}}_2+\bar{\omega}_3\bar{\mathcal{F}}_3|,\quad
\bar{\omega}'=(\bar{\omega}_2,\bar{\omega}_3)=\omega'/\omega,
\label{KK2}
\\
\label{ell0}
& \ell_0= \bar{\omega}_2\ell_2+\bar{\omega}_3\ell_3=\widebar{\mathcal{F}}_1\cdot(\bar{\omega}_2\widebar{\mathcal{F}}_2+
 \bar{\omega}_3\widebar{\mathcal{F}}_3).
\end{align}

Omitting straightforward calculations, we find the vector $\widetilde{U}_0=(\tilde{u}_1,\ldots ,\tilde{u}_{13})^{\top}$ in \eqref{3.2}, which is determined up to a nonzero factor:
\[
\tilde{u}_1=s-\frac{i\ell_0}{M^2}\omega ,\quad \tilde{u}_2=-d_0s+\frac{i\ell_0}{M^2}\omega,\quad
(\tilde{u}_3,\tilde{u}_4)=ia_0\omega',
\]
\[
(\tilde{u}_5,\tilde{u}_8,\tilde{u}_{11})^{\top}= -\Big(s-\frac{i\ell_0}{M^2}\omega\Big)\widebar{\mathcal{F}}_1 +\frac{ia_0\omega}{R}
(\bar{\omega}_2\bar{\mathcal{F}}_2+\bar{\omega}_3\bar{\mathcal{F}}_3),
\]
\[
(\tilde{u}_6,\tilde{u}_9,\tilde{u}_{12})^{\top}= ia_0\omega_2\widebar{\mathcal{F}}_1,\quad
(\tilde{u}_7,\tilde{u}_{10},\tilde{u}_{13})^{\top}= ia_0\omega_3\widebar{\mathcal{F}}_1.
\]
Then, we calculate the vector
\[
\mathcal{A}_1\widetilde{U}_0=-\frac{\beta}{2M}
\big(s\,,\;i\ell_0\omega-sM^2\,,\;iR(M^2-M_1^2)\omega'\,,\tilde{a}_5,0,0,
\tilde{a}_8,0,0,\tilde{a}_{11},0,0\big)^{\top}
\]
appearing in \eqref{3.6}, where $(\tilde{a}_5,\tilde{a}_8,\tilde{a}_{11})^{\top}=
-s\widebar{\mathcal{F}}_1+i\omega (\bar{\omega}_2\bar{\mathcal{F}}_2+\bar{\omega}_3\bar{\mathcal{F}}_3)$.

By replacing the first line of the matrix $s\mathcal{A}_0+\lambda^+ \mathcal{A}_1+{i}\omega_2 \mathcal{A}_2+{i}\omega_3 \mathcal{A}_3$ with the row-vector $(\mathcal{A}_1\widetilde{U}_0)^{\top}$, we get the matrix $\mathfrak{L}$ in \eqref{X}. We then compute the Lopatinski determinant
\[
\det \mathfrak{L} =\frac{\beta^2\Omega^6(\omega^2-\lambda^2)^2}{2M^2}\Big\{
(\lambda^2-\omega^2)s  +M^2s\Omega\lambda -2i\ell_0\omega\lambda^2  +M_{1}^2\lambda^2s
+ M_{2}^2\omega^2\lambda +R(M^2-M_{1}^2)\omega^2\Omega \Big\},
\]
where
$\lambda =\lambda^+$. Since
$\Omega^6(\omega^2-\lambda^2)^2\neq 0$ for $\lambda =\lambda^+$, the equality $\det \mathfrak{L}=0$ implies
\[
(\lambda^2-\omega^2)s  +M^2s\Omega\lambda -2i\ell_0\omega\lambda^2  +M_{1}^2\lambda^2s
+ M_{2}^2\omega^2\lambda +R(M^2-M_{1}^2)\omega^2\Omega =0,
\]
which can be rewritten as
\begin{equation}\label{84}
\Omega (M^2\lambda s +K\omega^2) +(M_*^2\lambda^2  -K_0\omega^2)s
=2i\ell_0\lambda^2\omega ,
\end{equation}
where
\begin{equation}\label{bigK}
K=R(M^2-M_1^2)+M_2^2
\end{equation}
($K>0$, cf. \eqref{Mach}).

The rest arguments are a suitable adaptation of the 2D spectral analysis in \cite{T22}. Following \cite{T22}, from the effective dispersion relation \eqref{83'} we obtain
\[
M_*^2\lambda^2 -K_0\omega^2=M^2\Omega^2-2i\ell_0\lambda\omega .
\]
By substituting the left-hand side of the last equality into \eqref{84}, one gets
\begin{equation}
\Omega (M^2\Omega^2-M^2\lambda^2+K\omega^2-2i\ell_0\lambda\omega )=0.
\label{0.}
\end{equation}
Since $\Omega\neq 0$ for $\lambda =\lambda^+$ and $\eta >0$, \eqref{0.} is reduced to
\begin{equation}
M^2\Omega^2-M^2\lambda^2+K\omega^2=2i\ell_0\lambda\omega .
\label{86=}
\end{equation}

Our analysis is thus reduced to the study of system \eqref{83'}, \eqref{86=}. As in \cite{T22}, without loss of generality we may suppose that $\ell_0 \geq 0$ (see \eqref{ell0}), i.e., we can finally consider the system
\begin{align}
 M^2\Omega^2-M_*^2\lambda^2 +K_0\omega^2 & =2i|\ell_0|\lambda\omega, \label{83}\\
 M^2\Omega^2-M^2\lambda^2+K\omega^2 & =2i|\ell_0|\lambda\omega .\label{86}
\end{align}
Indeed, if $\ell_0< 0$, then in \eqref{83'}, \eqref{86=} we make the change $\widetilde{\omega}'=-\omega' \in \mathbb{R}^2$.  After dropping tildes we again get system \eqref{83}, \eqref{86}. Subtracting \eqref{83} from \eqref{86} gives
\begin{equation}\label{l2}
\lambda^2=\frac{(K_0-K)\omega^2}{\beta^2}.
\end{equation}
This implies $\Real\lambda =0$ for $K\geq K_0$. Since $\Real\lambda^+ >0$ for $\eta >0$, it follows from \eqref{l2} that $\eta =0$ (for $K\geq K_0$). That is, we have proved

\begin{lemma}
In the parameter domain
\[
K\geq K_0
\]
the KL condition holds, i.e., shock waves cannot be violently unstable and they are, at least, weakly stable.
\label{lem1}
\end{lemma}

The effective dispersion relation \eqref{83} has the two roots
\[
\lambda_k=\frac{1}{\beta^2}\left( M^2s-i|\ell_0|\omega +(-1)^k\sqrt{M^2M_*^2s^2-2i|\ell_0|M^2s\omega+(K_0\beta^2-\ell_0^2)\omega^2} \right),\quad k=1,2.
\]
The Lax conditions \eqref{Mach} imply
\[
\Real\lambda_1|_{\omega' =0} =\frac{M}{\beta^2}(M-M_*)\eta <0,\quad \Real\lambda_2|_{\omega' =0} =\frac{M}{\beta^2}(M+M_*)\eta >0\quad\mbox{for}\quad \eta>0.
\]
It follows from Hersh's lemma \cite{Hersh} that $\Real\lambda_2 >0$ for all $\omega'\in\mathbb{R}^2$, i.e.,
\begin{equation}
\lambda^\pm=\frac{1}{\beta^2}\left( M^2s-i|\ell_0|\omega \pm\sqrt{M^2M_*^2s^2-2i|\ell_0|M^2s\omega+(K_0\beta^2-\ell_0^2)\omega^2}\right),
\label{l+}
\end{equation}
where $\lambda^-:=\lambda_1$.

For $\omega =0$, \eqref{l2} and \eqref{l+} imply $s =0$. Therefore, the UKL condition holds for the 1D case because $|s|^2+\omega^2 \neq 0$ (see Definition \ref{d3.3}; recall that $\omega:=|\omega'|$). This why we may assume that $\omega \neq 0$. Since the left-hand sides in \eqref{83} and \eqref{86} are homogeneous functions of $s$, $\lambda$ $\omega_2$ and $\omega_3$, without loss of generality, from now on we will suppose that $\omega =1$.

\begin{lemma}
In the parameter domain
\[
K< K_0
\]
the UKL condition holds, i.e., shock waves are uniformly stable.
\label{lem2}
\end{lemma}

The proof of Lemma \ref{lem2} follows from the fact that for $K< K_0$  system \eqref{83}, \eqref{86} necessarily imply $\eta <0$ (with $\lambda^+=\sqrt{K_0-K}/\beta\in\mathbb{R}$, cf. \eqref{l2}). Since the proof is really technical, we put it in Appendix \ref{Ap A}.

That is,  as follows from Lemmas \ref{lem1} and \ref{lem2}, the inequality
$K< K_0$ describes at least a part of the whole parameter domain of uniform stability. For finding the boundary of this domain we, roughly speaking, move along the $K$-axis from the point $K=K_0$ to the right and catch the point of transition to weak stability. The arguments are essentially similar to those for the 2D case in \cite{T22}.

Passing to the limit $\eta\rightarrow +0$ in \eqref{l+} (with $\omega =1$), we find
\begin{equation}\label{iml}
\Imag\lambda^\pm =\delta^\pm =\frac{1}{\beta^2}\left( M^2\xi-|\ell_0| \pm\sgn (M_*^2\xi -|\ell_0| ) \sqrt{M^2M_*^2\xi^2-2|\ell_0|M^2\xi+\ell_0^2-K_0\beta^2}\right).
\end{equation}
In the domain of uniform stability, system \eqref{83}, \eqref{86} must have no root $(s,\lambda )= (i\xi ,i\delta^+)$. Because of the continuous dependence of this system on $s$ and $\lambda$, the passage to weak stability (when the KL condition holds but the UKL condition is violated) may happen only thanks to the merging of $\delta^+$ and $\delta^-$ at a point of transition $\xi =\xi^*$. Clearly, for $\xi=\xi^*$ the square root in \eqref{iml} should vanish:
\begin{equation}
M^2M_*^2(\xi^*)^2-2|\ell_0|M^2\xi^*+\ell_0^2-K_0\beta^2 =0.
\label{xi2}
\end{equation}

Solving \eqref{xi2}, we find the two roots
\begin{equation}\label{xi*}
\xi^*=\xi^*_k :=\frac{M|\ell_0|+(-1)^{k+1}\beta\sigma}{MM_*^2},\quad k=1,2,
\end{equation}
which gives
\begin{equation}\label{d*}
\delta^+_{|\xi =\xi^*_k}=\delta^-_{|\xi =\xi^*_k}=  \delta^*_k:=\frac{M^2\xi^*_k-|\ell_0|}{\beta^2}\qquad (\mbox{for}\ k=1,2),
\end{equation}
where
\[
\sigma = \sqrt{M_*^2(1+M_2^2)-\ell_0^2}=\sqrt{1+|\widebar{\mathcal{F}}_1|^2+
|\bar{\omega}_2\widebar{\mathcal{F}}_2+\bar{\omega}_3\widebar{\mathcal{F}}_3|^2+
\big|\widebar{\mathcal{F}}_1\times (\bar{\omega}_2\widebar{\mathcal{F}}_2+\bar{\omega}_3\widebar{\mathcal{F}}_3)\big|^2}.
\]
It follows from \eqref{xi*} and \eqref{d*} that
\begin{equation}
\delta^*_1 =\frac{\sqrt{K_1}}{\beta}, \quad
\delta^*_2 =-\frac{\sqrt{K_2}}{\beta},
\label{delt}
\end{equation}
with
\[
K_1=\frac{(M\sigma -|\ell_0|\beta )^2}{M_*^4}>0,\quad
K_2=\frac{(M\sigma +|\ell_0|\beta )^2}{M_*^4}>K_1>0.
\]
Note that, in view of \eqref{Mach},  $K_1$ cannot vanish because
\begin{equation}
\begin{split}
M^2\sigma^2-\ell_0^2\beta^2=M_*^2(M^2+M^2M_2^2-\ell_0^2) & >M_*^2(M^2+M_1^2M_2^2-\ell_0^2) \\ & =
M_*^2\big(M^2+\big|\widebar{\mathcal{F}}_1\times (\bar{\omega}_2\widebar{\mathcal{F}}_2+\bar{\omega}_3\widebar{\mathcal{F}}_3)\big|^2\big)>0.
\end{split}
\label{K1}
\end{equation}
For $K>K_0$ (cf. Lemma \ref{lem2}), from \eqref{l2} and \eqref{delt} we find two possible transitions to weak stability:
\[
K=K_0+K_1\quad \mbox{and}\quad K=K_0+K_2
\]
($K_0+K_1<K_0+K_2$).

Since $\delta^\pm$ is real, the elementary analysis of the quadratic function in \eqref{xi2} shows that
\[
\xi \geq \xi_1^*\qquad\mbox{or}\qquad \xi\leq\xi_2^*.
\]
One can easily show that
\begin{equation}
\sgn (M_*^2\xi -|\ell_0| ) =\left\{\hspace*{-2mm}\begin{array}{r}
1\quad\mbox{for }\xi \geq \xi_1^*,\\[6pt]
-1\quad\mbox{for }\xi \leq \xi_2^*.
\end{array}
\right.
\label{sgn}
\end{equation}
In the domain of weak stability system \eqref{83}, \eqref{86} necessarily has the root  $(s,\lambda )= (i\xi ,i\delta^+)$. Taking \eqref{sgn} into account, we have
\[
\delta^+ \geq \frac{M^2\xi-|\ell_0|}{\beta^2}\geq \frac{M^2\xi_1^* -|\ell_0|}{\beta^2}=\delta^*_1>0 \qquad\mbox{for }\xi \geq \xi_1^*
\]
and
\[
\delta^+ \leq \frac{M^2\xi-|\ell_0|}{\beta^2}\leq \frac{M^2\xi_2^* -|\ell_0|}{\beta^2}=\delta^*_2<0 \qquad\mbox{for }\xi \leq \xi_2^*.
\]
implying that
\begin{align}
(\delta^+)^2\geq \frac{K_1}{\beta^2} &\qquad\mbox{for }\xi \geq \xi_1^*,
\label{KK1}
\\
(\delta^+)^2\geq \frac{K_2}{\beta^2} &\qquad\mbox{for }\xi \leq \xi_2^*.
\label{KK3}
\end{align}

In view of \eqref{l2}, the inequalities for $\delta^+$ in \eqref{KK1} and \eqref{KK3} are rewritten as
\[
K\geq K_0+K_1\quad\mbox{and}\quad K\geq K_0+K_2
\]
respectively. That is, moving along the $K$-axis from the point $K=K_0$ to the right, we meet the first point $K=K_0+K_1$ of transition from uniform to weak stability. Since for $K\geq K_0+K_2$ we still have has the root  $(s,\lambda )= (i\xi ,i\delta^+)$, the point $K= K_0+K_2$ is a fictitious point of transition, which just corresponds to the merging of $\delta^+$ and $\delta^-$. Hence, the parameter domain
\begin{equation}\label{KK}
K< K_0+K_1
\end{equation}
is that of uniform stability, whereas for $K\geq K_0+K_1$ shock waves are weakly stable. Inequality \eqref{KK}, which should be satisfied for all points on the circle $\bar{\omega}_2^2+\bar{\omega}_3^2=1$, is equivalently rewritten as \eqref{usc'}. This completes the proof of Theorem \ref{t1}.

\section{Concluding remarks}
\label{s6}

As in the 2D case in \cite{MTT20}, we can show that all compressive shock waves with convex equations of state $p=p(\rho )$ are uniformly stable, i.e., they are Lax shocks and inequality \eqref{usc'} is automatically fulfilled for them. Indeed, the proof of the following proposition is totally the same as in the 2D case in \cite{MTT20} and therefore we just omit it here.

\begin{proposition}
Let $p(\rho )$ be a convex function of $\rho$. Then all planar compressive shock waves satisfy the Lax conditions \eqref{Mach-} and \eqref{Mach} as well as  the  ``elastic'' counterpart
\begin{equation}
\widetilde{M}^2(R-1)<1
\label{gas'}
\end{equation}
of Majda's condition \eqref{gas}, where
\[
\widetilde{M}=\sqrt{M^2-M_1^2}=\sqrt{M^2-|\widebar{\mathcal{F}}_1|^2}
\]
is the ``elastic'' Mach number ($0<\widetilde{M}<1$, cf. \eqref{Mach}).
\label{pr1}
\end{proposition}

Omitting simple algebra, we equivalently rewrite the uniform stability condition \eqref{usc'}/\eqref{KK} as
\begin{equation}\label{usc1}
\widetilde{M}^2(R-1)<1+ \frac{\mathcal{D}}{M_*^4},
\end{equation}
where
\begin{equation}\label{DD1}
\mathcal{D}=(M\sigma -|\ell_0|\beta -M_*^2\widetilde{M})(M\sigma -|\ell_0|\beta +M_*^2\widetilde{M}).
\end{equation}
The proof of the following proposition is essentially the same as the proof of the corresponding inequality in the 2D case (see \cite[Appendix B]{MTT20}), but since there are still technical peculiarities of the 3D case, we present this proof in Appendix \ref{Ap B} for completeness.

\begin{proposition}
All planar Lax shock waves satisfy the inequality
\begin{equation}\label{DD}
\mathcal{D}>0.
\end{equation}
\label{pr2}
\end{proposition}

As in isentropic gas dynamics, rarefaction shocks ($R<1$) are, in general, possible. In view of \eqref{DD}, the uniform stability condition \eqref{usc1} holds for them. On the other hand, \eqref{usc1} is rewritten as
\begin{equation}\label{usc2}
M^2(R-1) < 1 + Q,
\end{equation}
where $Q=M_1^2 (R-1)+\frac{\mathcal{D}}{M_*^4}$. Thanks to \eqref{DD} the ``elastic'' additive $Q>0$ for $R>1$. Therefore, comparing \eqref{usc2} with Majda's condition \eqref{gas}, we make the conclusion that the elastic force plays {\it stabilizing role} for uniform stability.

As in the 2D case (see \cite[Theorem 2]{MTT20}), we come to the following theorem, which is the consequence of Propositions \ref{pr1} and \ref{pr2}.

\begin{theorem}
All compressive shock waves in isentropic elastodynamics with a convex equation of state $p=p(\rho )$ are structurally stable.
\label{t3}
\end{theorem}

It is natural to expect that the elastic force plays stabilizing role not only for uniform stability but also for the possible structural stability of weakly stable shocks. Since, as was shown by Coulombel  and Secchi \cite{CS}, weakly stable shock waves in isentropic gas dynamics are structurally stable, the natural conjecture is that weakly stable shock waves in isentropic elastodynamics are also structurally stable. The proof of this conjecture is postponed to the future research.

Regarding the sufficient structural stability condition \eqref{strc} (for general equations of state $p=p(\rho )$), its proof follows from the form of the variable coefficient linearized problem for nonplanar shocks. Indeed, for the system of conservation laws \eqref{7},
\[
\partial_tf^0(U)+\sum_{j=1}^{3}\partial_jf^j(U )=0,
\]
with $f^{\alpha}=f^{\alpha}(U)=(f_1^{\alpha}, \ldots ,f_{13}^{\alpha})^{\top}$ ($\alpha =\overline{0,3}$),
the Rankine-Hugoniot conditions \eqref{RH} read
\begin{equation}
B(U^+,U^-)\Phi -[f^1(U)]=0\quad \mbox{on }\Gamma ,
\label{a1}
\end{equation}
where $\Phi :=(\p_t\varphi,\p_2\varphi,\p_3\varphi)^{\top}$ and the matrix $B=B(U^+,U^-)$ of order $13\times 3$ is determined from the relation
\[
B(U^+,U^-)\Phi =\p_t\varphi [f^0(U)]+\p_2\varphi [f^2(U)]+\p_3\varphi [f^3(U)].
\]
After straightening  $\Gamma$, the linearization of \eqref{a1} about a basic state $(\widehat{U}^+(t,x),\widehat{U}^-(t,x),\hat{\varphi} (t,x'))$ results in the
variable coefficient boundary conditions
\begin{equation}
B(\widehat{U}^+,\widehat{U}^-)\Phi -[\widetilde{A}_{1}(\widehat{U},\hat{\varphi} )U]=0\quad \mbox{on }x_1=0,
\label{a2}
\end{equation}
where $\Phi$ is the same as above, but $\varphi$ is now the perturbation on $\hat{\varphi}$, the boundary matrix $\widetilde{A}_{1}(U,\varphi )$ is given in \eqref{AN},
\[
[\widetilde{A}_{1}(\widehat{U},\hat{\varphi} )U]:=\widetilde{A}_{1}(\widehat{U}^+,\hat{\varphi}^+ )U^+ - \widetilde{A}_{1}(\widehat{U}^-,\hat{\varphi}^- )U^-\quad \mbox{on }x_1=0,
\]
and $U^\pm$ is the perturbation of $\widehat{U}^\pm$. It is the boundary conditions in form \eqref{a2} (with a source term in the right-hand side) that are used in \cite{Tsiam} in the proof of the existence of solutions to the original nonlinear free boundary problem for nonplanar shocks by a fixed-point argument.

Freezing the coefficients in \eqref{a2} and setting $U^-\equiv 0$ for 1-shocks, we obtain the frozen coefficient linear boundary conditions
\begin{equation}
B(\widehat{U}^+,\widehat{U}^-)\Phi -\widetilde{A}_{1}(\widehat{U}^+,\hat{\varphi} )U=0\quad \mbox{on }x_1=0
\label{a3}
\end{equation}
for the frozen coefficient linear system
\begin{equation}
A_0(\widehat{U}^+)\partial_tU+\widetilde{A}_{1}(\widehat{U}^+,\hat{\varphi} )\partial_1U+A_2(\widehat{U}^+)\partial_2U+A_3(\widehat{U}^+)\partial_3U=0\qquad \mbox{for}\ x\in \mathbb{R}^3_{\pm}
\label{23'}
\end{equation}
associated with \eqref{23}, where $U:=U^+$. Taking  into account the form of the boundary matrix $\widetilde{A}_{1}(U,\varphi )$, after the reduction of problem \eqref{a3}, \eqref{23'} to a dimensionless form we get a linear problem in form \eqref{ls}/\eqref{lsm}, \eqref{lbc}  in which $M$, $\mathcal{F}_{1i}$ and $\mathcal{F}_{kj}$ ($i,j=1,2,3$, $k=2,3$) are replaced with $\mathcal{M}$, $F^+_{iN}/(c^+|N|)$ and $F^+_{kj}/c^+$  respectively (see \eqref{ff}) written on the frozen basic state $(\widehat{U}^+,\widehat{U}^-,\hat{\varphi})$. It is clear that the uniform stability condition for this linear constant coefficient problem is nothing else but \eqref{strc} following from \eqref{usc'} after suitable replacements. This implies Theorem \ref{t2}.

At last, we note that structural stability in the sense of local-in-time existence in Theorem \ref{t2} can be alternatively replaced with structural stability in the sense that the nonlinear free boundary problem for a shock wave being {\it initially closed to a planar shock} has a unique smooth solution on a finite (but not necessarily small) time interval (see, e.g., \cite{CS}).

\section*{Acknowledgements}

The author Y. Trakhinin acknowledges the financial support from the Russian Science Foundation (Project No. 24-21-00192).

\appendix

\section{Proof of Lemma 5.2}

\label{Ap A}

For $K< K_0$,  it follows from \eqref{l2} (for $\omega =1$) that
\begin{equation}
\label{A.1}
\lambda^+=\frac{\sqrt{K_0-K}}{\beta}\in\mathbb{R}.
\end{equation}
Considering the imaginary and real  parts of $\eqref{86}$, we obtain
\begin{equation}
\label{A.2}
\xi = \frac{|\ell_0|\lambda^+}{M^2(\eta + \lambda^+)}
\end{equation}
and
\begin{equation}
\label{A.3}
M^2(\eta^2- \xi ^2 + 2\eta\lambda^+) + K = 0,
\end{equation}
respectively.
Substituting \eqref{A.2} into \eqref{A.3}, we get
\begin{equation}
\label{A.4}
M^4\gamma^4 - (M^4(\lambda^+)^2 - KM^2)\gamma^2 - |\ell_0|^2(\lambda^+)^2 = 0,
\end{equation}
where $\gamma = \eta + \lambda^+$. Solving the quadratic equation \eqref{A.4}, we get
\begin{equation}
\label{A.5}
\gamma^2 = \frac{1}{2M^2}\biggl(M^2(\lambda^+)^2 - K + \sqrt{\bigl(M^2(\lambda^+)^2 - K\bigr)^2 + 4|\ell_0|^2(\lambda^+)^2} \biggr).
\end{equation}

Let $\theta$ denote the right-hand side of \eqref{A.5}, then $\gamma^2 = \theta$ implies
\[
\eta = \pm \sqrt{\theta} - \lambda^+ < \sqrt{\theta} - \lambda^+.
\]

We now prove that $\eta <0$ (for $K< K_0$), i.e., $\sqrt{\theta} - \lambda^+ < 0$. By virtue of \eqref{A.1}, the last inequality is equivalent to $\theta < (\lambda^+)^2$ or
\begin{equation*}
\frac{1}{2M^2}\biggl(M^2(\lambda^+)^2 - K + \sqrt{\bigl(M^2(\lambda^+)^2 - K\bigr)^2 + 4|\ell_0|^2(\lambda^+)^2} \biggr) < (\lambda^+)^2.
\end{equation*}
The latter is equivalent to
\begin{equation}
\label{A.6}
KM^2-\ell_0^2>0.
\end{equation}

By virtue of \eqref{KK2}, \eqref{ell0} and \eqref{bigK}, inequality \eqref{A.6} is rewritten as
\begin{equation}
\label{A.7}
R(M^2-M_1^2)M^2+ |\bar{\omega}_2\bar{\mathcal{F}}_2+\bar{\omega}_3\bar{\mathcal{F}}_3|^2M^2-
\big(\widebar{\mathcal{F}}_1\cdot(\bar{\omega}_2\widebar{\mathcal{F}}_2+
 \bar{\omega}_3\widebar{\mathcal{F}}_3)\big)^2>0.
\end{equation}
In view of the Lax condition $M^2>M_1^2=|\widebar{\mathcal{F}}_{1}|^2$ (see \eqref{Mach}), inequality \eqref{A.7} holds as soon as
\[
|\widebar{\mathcal{F}}_{1}|^2|\bar{\omega}_2\bar{\mathcal{F}}_2+\bar{\omega}_3\bar{\mathcal{F}}_3|^2-
\big(\widebar{\mathcal{F}}_1\cdot(\bar{\omega}_2\widebar{\mathcal{F}}_2+
 \bar{\omega}_3\widebar{\mathcal{F}}_3)\big)^2=\big|\widebar{\mathcal{F}}_1\times (\bar{\omega}_2\widebar{\mathcal{F}}_2+\bar{\omega}_3\widebar{\mathcal{F}}_3)\big|^2>0.
\]
Since the letter is always true, we conclude that the Lopatinski determinant vanishes only for $\eta <0$. This implies that the UKL condition holds for $K<K_0$.

\section{Proof of Propisition 6.2}

\label{Ap B}

The proof of Proposition \ref{pr1} is just a suitable ``3D correction'' of the proof of the 2D analogue of inequality \eqref{DD} in \cite[Appendix B]{MTT20}.

By virtue of \eqref{K1} and \eqref{DD1}, the condition  \eqref{DD} is equivalent to
\[
M\sigma -|\ell_0|\beta - M_*^2\widetilde{M}>0.
\]
or
\begin{equation}
M\sigma - M_*^2\widetilde{M}>|\ell_0|\beta .
\label{D'}
\end{equation}
The left-hand side in \eqref{D'} is strictly positive because
\[
M^2\sigma^2 -M_*^4\widetilde{M}^2=M_*^2M_1^2\beta^2 +(M_2^2+\varkappa^2 )M^2>0,
\]
where $\varkappa =\big|\widebar{\mathcal{F}}_1\times (\bar{\omega}_2\widebar{\mathcal{F}}_2+\bar{\omega}_3\widebar{\mathcal{F}}_3)\big| >0$ (the 2D counterpart of $\varkappa$ in \cite[Appendix B]{MTT20} is $\varkappa =\det \mathcal{F}$). Therefore,  squaring \eqref{D'} gives the equivalent inequality
\begin{equation}
M^2\sigma^2 -\ell_0^2\beta^2+M_*^4\widetilde{M}^2> 2M_*^2\sigma M\widetilde{M}.
\label{D''}
\end{equation}

In view of \eqref{K1}, the left-hand side in \eqref{D''} is strictly positive. Squaring \eqref{D''}, we again obtain an equivalent inequality which after long but straightforward calculations can be written as
\[
4\varkappa^2\widetilde{M}^2\big(1-\widetilde{M}^2\big) +\big( (M_1^2+M_2^2)\widetilde{M}^2-(M_1^2+\varkappa^2)\big)^2 +4M_2^2\varkappa^2\widetilde{M}^2>0.
\]
The last inequality is true because $0<\widetilde{M}<1$ (cf. \eqref{Mach}) for Lax shocks. We have thus proved that $\mathcal{D}>0$.

{\footnotesize 
  }

\end{document}